\documentclass[11pt]{article}

\usepackage{amsfonts,amsmath,latexsym}
\newtheorem{theorem}{Theorem}
\newtheorem{lemma}{Lemma}

\newenvironment{proof}
      {\medskip\noindent{\bf Proof:}\hspace{1mm}}
      {\hfill$\Box$\medskip}

\def\qed{\ifvmode\mbox{ }\else\unskip\fi\hskip 1em plus 10fill$\Box$}
\makeatletter
\def\Ddots{\mathinner{\mkern1mu\raise\p@
\vbox{\kern7\p@\hbox{.}}\mkern2mu
\raise4\p@\hbox{.}\mkern2mu\raise7\p@\hbox{.}\mkern1mu}}
\makeatother

\title{\vspace{-0.7cm}Paths and stability number in digraphs}
\author{Jacob
Fox\thanks{Department of Mathematics, Princeton, Princeton, NJ.
Email: {\tt jacobfox@math.princeton.edu}. Research supported by an
NSF Graduate Research Fellowship and a Princeton Centennial
Fellowship.} \and Benny Sudakov\thanks{Department of Mathematics,
UCLA,  Los Angeles, CA 90095. Email: {\tt bsudakov@math.ucla.edu}. Research supported in part by NSF CAREER award DMS-0812005 and by
USA-Israeli BSF grant.}}
   \date{}
\begin{document}
\maketitle

\begin{abstract}

The Gallai-Milgram theorem says that the vertex set of any digraph with stability number $k$ can be partitioned into $k$ directed paths. In 1990, Hahn and Jackson conjectured that this theorem is best
possible in the following strong sense. For each positive integer $k$, there is a digraph $D$ with stability number $k$ such that deleting the vertices of any $k-1$
directed paths in $D$ leaves a digraph with
stability number $k$. In this note, we prove this conjecture.
 \end{abstract}

\section{Introduction}

The Gallai-Milgram theorem \cite{GaMi} states that the vertex set of any digraph with stability number $k$ can be
partitioned into $k$ directed paths. It generalizes Dilworth's theorem \cite{Di} that the size of a maximum antichain in
a partially ordered set is equal to the minimum number of chains needed to cover it. In 1990, Hahn and Jackson
\cite{HaJa} conjectured that this theorem is best possible in the following strong sense. For each positive integer $k$,
there is a digraph $D$ with stability number $k$ such that deleting the vertices of any $k-1$ directed paths in $D$
leaves a digraph with stability number $k$. Hahn and Jackson used known bounds on Ramsey numbers to verify their
conjecture for $k \leq 3$. Recently, Bondy, Buchwalder, and Mercier \cite{BoMe} used lexicographic products of graphs to
show that the conjecture holds if $k=2^a3^b$ with $a$ and $b$ nonnegative integers. In this short note we prove the
conjecture of Hahn and Jackson for all $k$.

\begin{theorem}\label{main}
For each positive integer $k$, there is a digraph $D$ with stability number $k$ such that deleting the vertices of any $k-1$ directed paths leaves a
digraph with stability number $k$.
\end{theorem}

To prove this theorem we will need some properties of random graphs.
As usual, the random graph $G(n,p)$ is a graph on $n$ labeled vertices in which each pair of vertices forms an edge randomly and independently with
probability $p=p(n)$.

\begin{lemma}\label{mainlemma}
For $k \geq 3$, the random graph $G=G(n,p)$ with $p=20n^{-2/k}$ and $n \geq 2^{15k^2}$ a multiple of $2k$ has the following properties. \\
(a) The expected number of cliques of size $k+1$ in $G$ is at most $20^{{k+1 \choose 2}}$. \\
(b) With probability more than $\frac{2}{3}$, every induced subgraph of $G$ with $\frac{n}{2k}$ vertices has a clique of size $k$.
\end{lemma}
\begin{proof}
(a) Each subset of $k+1$ vertices has probability $p^{{k+1 \choose 2}}$ of being a clique. By linearity of expectation, the
expected number of cliques of size $k+1$ is
$${n \choose k+1}p^{{k+1 \choose 2}} = {n \choose k+1}20^{{k+1 \choose 2}}n^{-k-1} \leq 20^{{k+1 \choose 2}}.$$
(b) Let $U$ be a set of $\frac{n}{2k}$ vertices of $G$. We first give an upper bound on the probability that $U$ has no clique of size $k$.
For each subset $S \subset U$ with $|S|=k$, Let $B_S$ be the event that $S$ forms a clique, and $X_S$ be the indicator random variable for $B_S$.
Since $k \geq 3$, by linearity of expectation, the expected number $\mu$ of cliques in $U$ of size $k$ is $$\mu=\mathbb{E}\left[\sum_{S} X_S\right]={\frac{n}{2k} \choose k}p^{{k \choose 2}} \geq
\frac{n^k}{2(2k)^kk!}20^{{k \choose 2}}n^{1-k} \geq 2n.$$
Let $\Delta=\sum \Pr[B_S \cap B_T],$ where the sum is over all ordered pairs $S,T$ with $|S \cap T| \geq 2$. We have
\begin{eqnarray*}
\Delta & = & \sum_{i=2}^{k-1}\sum_{|S \cap T|=i} \Pr[B_S \cap B_T]=\sum_{i=2}^{k-1}\sum_{|S \cap T|=i} p^{2{k \choose 2}-{i \choose 2}}=\sum_{i=2}^{k-1}{n \choose i}{n-i \choose k-i}{n-k \choose
k-i}p^{2{k \choose 2}-{i \choose 2}} \\ & \leq &  \sum_{i=2}^{k-1} n^{2k-i}p^{k(k-1)-{i \choose 2}} \leq 20^{k^2} \sum_{i=2}^{k-1} n^{2-i+i(i-1)/k} \leq k20^{k^2}n^{2/k}\,.
\end{eqnarray*}
Here we used the fact that $i(i-1)/k-i$ for $2 \leq i \leq k-1$ clearly achieves its maximum when $i=2$ or $i=k-1$.

Using that $k \geq 3$ and $n \geq 2^{15k^2}$, it is easy to check that $\Delta \leq n$. Hence, by Janson's inequality (see, e.g., Theorem 8.11 of \cite{AlSp}) we can bound the probability that $U$ does not contain a clique of size $k$ by
$\Pr\left[ \wedge_{S} \bar{B}_S\right] \leq e^{-\mu+\Delta/2} \leq e^{-n}$.
By the union bound, the probability that there is a set of $\frac{n}{2k}$ vertices of $G(n,p)$ which does not contain a clique of size $k$ is at most
${n \choose \frac{n}{2k}}e^{-n} \leq 2^ne^{-n} < 1/3$.
\end{proof}

The proof of Theorem \ref{main} combines the idea of Hahn and Jackson of partitioning a graph
into maximum stable sets and orienting the graph accordingly with Lemma \ref{mainlemma} on properties of random graphs.

\noindent
{\bf Proof of Theorem \ref{main}.}
Let $k \geq 3$ and $n \geq 2^{15k^2}$. By Markov's inequality and Lemma \ref{mainlemma}(a), the probability
that $G(n,p)$ with $p=20n^{-2/k}$ has at most $2 \cdot 20^{{k+1 \choose 2}}$ cliques of size $k+1$ is at least $1/2$.
Also, by Lemma \ref{mainlemma}(b), we have that with probability at least $2/3$ every set of $\frac{n}{2k}$ vertices of this random graph contains a clique
of size $k$. Hence, with positive probability (at least $1/6$) the random graph $G(n,p)$ has both properties.
This implies that there is a graph $G$ on $n$ vertices which contains at most $2 \cdot 20^{{k+1 \choose 2}}$ cliques of size $k+1$ and every
set of $\frac{n}{2k}$ vertices of $G$ contains a clique of size $k$.
Delete one vertex from each clique of size $k+1$ in $G$. The resulting graph $G'$ has at least $n-2 \cdot 20^{{k+1 \choose 2}} \geq 3n/4$
vertices and no cliques of size $k+1$. Next pull out vertex disjoint cliques of size $k$ from $G'$ until the remaining subgraph has no clique of size $k$, and let
$V_1,\ldots,V_t$ be the vertex sets of these disjoint cliques of size $k$. Since every induced subgraph of $G$ of size at least $\frac{n}{2k}$ contains a clique of size $k$,
then $|V_1 \cup \ldots \cup V_t| \geq \frac{3n}{4}-\frac{n}{2k} \geq \frac{n}{2}$.
Define the digraph $D$ on the vertex set $V_1 \cup \ldots \cup V_t$ as follows.
The edges of $D$ are the nonedges of $G$. In particular, all sets $V_i$ are stable sets in $D$. Moreover, all edges of $D$ between $V_i$ and $V_j$ with $i<j$ are oriented from $V_i$ to $V_j$. By construction, the stability number of $D$ is equal to the clique number of $G'$, namely $k$. Also any set of $\frac{n}{2k}$ vertices of $D$ contains a stable set of size $k$. Note that every directed path in $D$
has at most one vertex in each $V_i$. Hence, deleting any $k-1$ directed paths in $D$ leaves at least $|D|/k \geq \frac{n}{2k}$ remaining vertices.
These remaining vertices contain a stable set of size $k$, completing the proof. \qed

\noindent
{\bf Remark.}\,
Note that in order to prove Theorem \ref{main}, we only needed to find
a graph $G$ on $n$ vertices with no clique of size $k+1$
such that every set of $\frac{n}{2k}$ vertices of $G$ contains a clique of
size $k$. The existence of such graphs were first proved by
Erd\H{o}s and Rogers \cite{ErRo}, who more generally asked to estimate
the minimum $t$ for which there is a graph $G$ on $n$ vertices with no clique of size $s$ such
that every set of $t$ vertices of $G$ contains a clique of size $r$.
Since then a lot of work has been done on this question, see, e.g.,
\cite{Kr,AlKr,Su,DuRo}. Although most result for this problem used probabilistic arguments,
Alon and Krivelevich \cite{AlKr} give an explicit construction
of an $n$-vertex graph $G$ with no clique of size
$k+1$, such that every subset of $G$ of size $n^{1-\epsilon_k}$ contains a
$k$-clique. Since we only need a much weaker result to prove the conjecture of Hahn and Jackson, we decided to include its very short and simple proof to keep this note self-contained.

\vspace{0.1cm} \noindent {\bf Acknowledgments.}\, We would like to thank Adrian Bondy for stimulating discussions and
generously sharing his presentation slides. We also are  grateful to Noga Alon for drawing our attention to the paper
\cite{AlKr}. Finally, we want to thank the referee for helpful comments.

\end{document}